\documentclass{amsart}
\usepackage{amssymb}


\def\R{{\mathbb R}}
\def\dach{^\wedge} 
\def\speccocx{(C_0(X)\rtimes G)\dach}
\def\restr#1{|_{{#1}}}
\def\set#1{\{\,#1\,\}} 
\def\K{{\mathcal K}}
\def\pedid{\kappa} 
\def\H{{\mathcal H}} 
\def\M{{m}} 
\def\TT{\mathcal{T}}
\def\O{{\mathcal O}}
 \def\calS{{\mathcal S}}
\def\hG{\widehat G}
\def\cocx{C_{0}(X)\rtimes G}
\def\cocy{C_{0}(Y)\rtimes G}

\makeatletter
\def\labelenumi{\textnormal{(\@alph\c@enumi)}}
\def\theenumi{\@alph \c@enumi}
\newcount\charno
\def\alphapart#1{\charno=96
\advance\charno by#1\char\charno}

\makeatother

\def\Ind{\operatorname{Ind}} 
\def\interior{\operatorname{int}}

\def\Prim{\operatorname{Prim}} 
\def\tr{\operatorname{tr}}

\newtheorem{thm}{Theorem}[section] 
\newtheorem{cor}[thm]{Corollary}
\newtheorem{lemma}[thm]{Lemma} 
\newtheorem{prop}[thm]{Proposition}

\theoremstyle{definition}

\theoremstyle{remark} 
\newtheorem{remark}[thm]{Remark}
\newtheorem{example}[thm]{Example}

\numberwithin{equation}{section}

\begin{document}

\title[Ideals in transformation-group $C^*$-algebras]{\boldmath Ideals
  in transformation-group $C^*$-algebras}

\author[an Huef]{Astrid an Huef} 

\author[Williams]{Dana P. Williams}

\address{Department of Mathematics and Computer Science\\
  University of Denver\\
  2360 S. Gaylord St.\\
  Denver, CO 80208\\
  USA} \email{astrid@cs.du.edu}

\address{Department of Mathematics\\
  Dartmouth College\\
  Hanover, NH 03755\\
  USA} \email{dana.williams@dartmouth.edu}

\thanks{The authors thank Judy Packer for a very helpful
  discussion.}

\date{June 5, 2000}

\begin{abstract} 
  We characterize the ideal of continuous-trace elements in a
  separable transformation-group $C^{*}$-algebra $\cocx$.  In
  addition, we identify the largest Fell ideal, the largest liminal
  ideal and the largest postliminal ideal.
\end{abstract}

\maketitle

 \section{Introduction}
 
 Let $(G,X)$ be a locally compact Hausdorff transformation group: thus
 $G$ is a locally compact Hausdorff group and $X$ is a locally compact
 Hausdorff space together with a jointly continuous map $(s,x)\mapsto
 s\cdot x$ from $G\times X$ to $X$ such that $s\cdot (t\cdot
 x)=st\cdot x$ and $e\cdot x=x$.  The associated transformation-group
 $C^*$-algebra $C_0(X)\rtimes G$ is the $C^*$-algebra which is
 universal for the covariant representations of the $C^*$-dynamical
 system $(C_0(X), G, \alpha)$ in the sense of \cite{rae88:_takai}.
 More concretely, $C_0(X)\rtimes G$ is the enveloping $C^*$-algebra of
 the Banach $*$-algebra $L^1(G, C_0(X))$ of functions $f:G\to C_0(X)$
 which are Bochner integrable with respect to a fixed left Haar
 measure on $G$ \cite[Section~7.6]{ped79:_automorphisms}. We
 will always assume that $G$ and $X$ are second countable so that
 $\cocx$ is separable.  In our main results, we assume either that $G$
 is abelian or that $G$ acts freely.

 It is natural to attempt to characterize properties of $\cocx$ in
 terms of the dynamics of the action of $G$ on $X$, and there are a
 large number of results of this sort in the literature
 \cite{goot:_type,green77:_orbit,
  will81:_ccr,will81:_tr,will82:_spectrum,aah:trace,aah:fell}.
 We were motivated by a particularly nice example due to Green
 \cite[Corollary~18]{green77:_orbit} 
 in which he was able to characterize the closure $I$ of the ideal of
 continuous-trace elements in $\cocx$ in the case $G$ acts freely and
 $\cocx$ is postliminal.  (Since we'll be working exclusively with
 separable $C^{*}$-algebras, we will not distinguish between Type~I
 and postliminal algebras.)  There are three ingredients required for
 this sort of project.  First, one needs a global characterization of
 algebras $\cocx$ which have continuous trace.  Second, one needs to
 know that the ideal $I$ is of the form $C_{0}(Y)\rtimes G$ for an
 open $G$-invariant set $Y$ in $X$.  And third, one wants a
 straightforward characterization of $Y$ in terms of the dynamics.
 Assuming that $G$ acts freely, Green showed that $\cocx$ has
 continuous trace if and only if every compact set $K\subset X$ is
 \emph{wandering} in that
\begin{equation*}
  \set{s\in G: s\cdot K \cap K\not=\emptyset}
\end{equation*}
is relatively compact in $G$ \cite[Theorem~17]{green77:_orbit}.  If
$\cocx$ is postliminal, then $G$ acting freely implies every ideal is
of the form $C_{0}(Y)\rtimes G$, and Green showed $I=C_{0}(Y)\rtimes G$
where
\begin{multline}\label{eq:7}
  Y=\{\,y\in X: \text{$y$ has a compact wandering neighborhood $N$} \\
  \text{such that $q(N)$ is closed and Hausdorff}\,\},
\end{multline}
where $q:X\to X/G$ is the quotient map.  (The criteria in \eqref{eq:7}
are slightly different than those given by Green; unfortunately, the
statement in \cite[Corollary~18]{green77:_orbit} is not quite correct
--- see Remark~\ref{rem-fu}.)

To extend Green's results to actions which are not necessarily free,
we relied (1)~on the second author's result
\cite[Theorem~5.1]{will81:_tr} stating that if $G$ is abelian then
$\cocx$ has continuous trace if and only if the stability groups move
continuously and every compact set is \emph{$G$-wandering} as defined
in \S\ref{sec:ident-ideals}, and (2)~on a result of Phillips which
allows us to assume the ideal in question is of the form
$C_{0}(Y)\rtimes G$.  Our characterization is given in
Theorem~\ref{thm-new-main} and is valid for abelian groups, freely
acting amenable groups, or freely acting groups for which $\cocx$ is
postliminal.

For abelian groups or freely acting groups, Gootman showed that
$\cocx$ is postliminal if and only if the orbit space $X/G$ satisfies
the $\mathrm{T}_{0}$ axiom of separability
\cite[Theorem~3.3]{goot:_type}.  Similarly $\cocx$ is liminal if and
only if each orbit is closed \cite[Theorem~3.1]{will81:_ccr}.  Using
these results, we give characterizations of the largest postliminal
and liminal ideals in $\cocx$ in Theorems \ref{thm-I-ideal}~and
\ref{thm-liminal-ideal}, respectively.

The set of $a\in A^{+}$ such that $\pi\mapsto\tr\pi(a)$ is bounded on
$\hat A$ is the positive part of a two-sided ideal $\TT(A)$.  If
$\TT(A)$ is dense in $A$, then $A$ is said to have \emph{bounded
  trace}.  Such algebras are also \emph{uniformly liminal}
\cite[Theorem~2.6]{ass:97}.  The first author has characterized when
$\cocx$ has bounded trace \cite[Theorem~4.9]{aah:trace}, and she has
used this to find the largest bounded trace ideal in
\cite[Theorem~5.8]{aah:trace}.  An intermediate condition between $A$
being a continuous-trace $C^{*}$-algebra and an algebra with bounded
trace is that $A$ be a \emph{Fell algebra}.  A point $\pi\in\hat A$ is
called a \emph{Fell point} of the spectrum if there is a neighborhood
$V$ of $\pi$ and $a\in A^{+}$ such that $\rho(a)$ is a rank-one
projection for all $\rho\in V$.  Then $A$ is a Fell algebra if every
$\pi\in \hat A$ is a Fell point, and a Fell algebra is a
continuous-trace $C^{*}$-algebra if and only if $\hat A$ is Hausdorff
(cf., \cite[\S5.14]{rw-book}).  If $G$ acts freely, then $\cocx$ is a
Fell algebra if and only if $X$ is a Cartan $G$-space \cite{aah:fell},
and we treat the case of continuously varying stabilizers below
(Proposition~\ref{prop-fell}).  Using these results, we identify the
largest Fell ideal in $\cocx$ when the stability groups vary
continuously (Corollary~\ref{cor-fell}).

Naturally our techniques depend on describing ideals in $\cocx$ in
terms of the dynamics.  To do this, we need to know that each
primitive ideal in $\cocx$ is induced from a stability group (cf.,
\cite[Definition~4.12]{will81:_ccr}).  Cross products with this
property are called \emph{EH-regular}, and in the separable case it
suffices for $G$ to be amenable \cite{goot-ros} or for the orbit space
$X/G$ to be $\mathrm{T}_{0}$ \cite[Proposition~20]{green78}.
Therefore, if $G$ is abelian then $\cocx$ is EH-regular.
If $G$ acts freely, then we will have to assume either that $G$ is
amenable or the orbit space is $\mathrm{T}_{0}$. 
If the action
is free and $\cocx$ is EH-regular, then ideals in $\cocx$ are in
one-to-one correspondence with $G$-invariant open sets $Y$ in $X$.  If
$G$ does not act freely, then we must assume that $G$ is abelian so
that we can employ the dual action to conclude that the ideals we are
interested in correspond to $G$-invariant open subsets of $X$.

\vfill

\section{Invariance of ideals under the dual action}

Although ideals in $\cocx$ can be difficult to describe in general,
there is always an ideal associated to each $G$-invariant open subset
$Y$ of $X$.  The closure of $C_{c}(G\times Y)$ (viewed as a subset of
$C_{c}(G\times X)$) is an ideal in $\cocx$ which we can identify with
$C_{0}(Y)\rtimes G$ (cf., e.g., \cite[Lemma~1]{green77:_orbit}).  When
the action of $G$ is free and $\cocx$ is EH-regular,
\cite[Corollary~5.10]{will81:_ccr} implies that
$\Prim\bigl(\cocx\bigr)$ is homeomorphic to the quotient space
$(X/G)^{\sim}$ of $X/G$ where $G\cdot x$ is identified with $G\cdot y$
if $\overline{G\cdot x}=\overline{G\cdot y}$.  It follows that every
ideal of $\cocx$ is of the form $C_{0}(Y)\rtimes G$ for some
$G$-invariant open set $Y$.  When $G$ is abelian and does not
necessarily act freely, we can distinguish those ideals of $\cocx$ of
the form $C_{0}(Y)\rtimes G$ via the dual action.  Indeed, let $\hG$
denote the Pontryagin dual of $G$.  The \emph{dual action}
$\hat\alpha$ of $\hG$ on $C_0(X)\rtimes G$ is given by
\begin{equation*}
\hat\alpha_\tau(f)(s)=\tau(s)f(s,\cdot)\quad \text{for $f\in
C_c(G\times X)$ and $\tau\in\hG$}.
\end{equation*}
The induced action of $\hG$ on $(C_0(X)\rtimes G)\dach$ is
$\tau\cdot\pi=\pi\circ\hat\alpha_\tau^{-1}$, and this action is
jointly continuous (cf., e.g., \cite[Lemma~7.1]{rw-book}).  The
importance of the dual action for us comes from the following lemma
due to Phillips.

\begin{lemma} [{\cite[Proposition~6.39]{phillips88:_equivariant}}]
  \label{lem-invar-ideals} 
  Suppose that $(G,X)$ is a second countable transformation group with
  $G$ abelian.  If $I$ is a $\hG$-invariant ideal of $C_{0}(X)\rtimes
  G$, then there is an open $G$-invariant set $Y$ in $X$ such that $I=
  C_{0}(Y) \rtimes G$.
\end{lemma}

As an example, note that it is easy to see that the set of Fell points
of the spectrum is invariant under the dual action. If $\pi$ is a Fell
point, then by definition there exist $a\in A^+$ and an open
neighborhood $V$ of $\pi$ in $\hat A$ such that $\sigma(a)$ is a
rank-one projection for all $\sigma\in V$. If $b=\hat\alpha_\tau(a)$
then for every $\rho\in\tau\cdot V$ we have $\rho(b)=\sigma(a)$ for
some $\sigma\in V$. Hence $\tau\cdot\pi$ is also a Fell point.  Thus
the largest Fell ideal must be of the form $C_{0}(Y)\rtimes G$.

Recall that a positive element $a$ of a $C^*$-algebra $A$ is a
\emph{continuous-trace element} if the function
$\pi\mapsto\tr(\pi(a))$ is finite and continuous on $\hat A$. 
The linear span $\M(A)$ of these  elements is an ideal in $A$, and $A$ is a 
\emph{continuous-trace $C^*$-algebra} if  $\M(A)$ is  dense in $A$.

We want to prove that $\overline{\M(A)}$ is invariant under the dual
action.  To do this, we need a lemma of Green which characterizes this
ideal by determining its irreducible representations.  Recall that if
$I$ is an ideal of a $C^*$-algebra $A$, then the spectrum $\hat I$ of
$I$ is homeomorphic to the open set $\O_I:=\set{\rho\in\hat
  A:\rho(I)\not=\set0}$ in $\hat A$.  We will also use that every
$C^*$-algebra $A$ has a dense hereditary ideal $\pedid(A)$ --- called
the \emph{Pedersen ideal} of $A$ --- which is the smallest dense ideal
in $A$ \cite[Theorem~5.6.1]{ped79:_automorphisms}.  As Green's result
is an essential ingredient in many of our proofs, we give the brief
argument here.  The key idea of the proof is that
$\pi\bigl(\M(A)\bigr) \not=\set0$ if and only if $\pi$ has lots of
\emph{closed} neighborhoods in $\hat A$.

\begin{lemma}[{\cite[p.~96]{green77:_orbit}}] \label{lem-green}%
  Let $A$ be a $C^*$-algebra and $I=\overline{\M(A)}$. Then
  $\pi\in\O_I$ if and only if
  \begin{enumerate}
  \item there exists an ideal $J$ of $A$ which has continuous trace
    such that $\pi\in\O_J$\label{item:1}; and
  \item $\pi$ has a neighborhood basis consisting of closed
    sets.\label{item:2}
  \end{enumerate}
\end{lemma}

\begin{proof}
  Let $\pi\in\O_I$. There exists a positive element
  $a\in\M(A)$ such that $\tr(\pi(a))~=~1$. It follows that
  the set
\begin{equation*}
L=\set{\rho\in\hat A : \tr(\rho(a))\geq\textstyle{\frac{1}{2}}}
\end{equation*}
is a closed neighborhood of $\pi$ and $L\subset\O_I$.  Let $\set{
  F_\alpha}$ be a compact neighborhood basis of $\pi$ in $\hat A$.
Notice that $L$ is Hausdorff since $\mathcal{O}_{I}$ is.  Thus
$F_{\alpha}\cap L$ is closed in $L$, and therefore in $\hat A$ as
well.  It follows that $\set{ F_\alpha\cap L}$ is a neighborhood basis
of $\pi$ consisting of closed sets. This proves item~(\ref{item:2}).
That item~(\ref{item:1}) holds is obvious (just take $J=I$).

Conversely, let $\pi\in\hat A$ satisfy items~(\ref{item:1})
and~(\ref{item:2}).  Then there exists an ideal $J_0\subset J$ of $A$
such that $\pi\in\O_{J_0}$ and $\overline{\O_{J_0}}\subset \O_J$.  Let
$a$ be a positive element of the Pedersen ideal $\pedid(J_0)$ of
$J_0$. Then $\rho\mapsto\tr(\rho(a))$ is continuous on $\O_J$ because
$\pedid(J_0)\subset\pedid(J)\subset\M(J)$. Since
$\rho\mapsto\tr(\rho(a))$ vanishes off of $\O_{J_0}$, it is continuous
on all of $\hat A$. Thus $\pedid(J_0)\subset\M(A)\subset I$, whence
$J_0\subset I$, and $\pi\in\O_I$.
\end{proof}

\begin{prop}\label{prop-cts-trace-id}
  Let $(G,X)$ be a second countable transformation group with $G$
  abelian.  Then $I=\overline{\M(C_0(X)\rtimes G)}$ is
  $\hG$-invariant, and $I=C_0(Y)\rtimes G$ for an open $G$-invariant
  subset $Y$ of $X$.
\end{prop}

\begin{proof}We use Lemma~\ref{lem-green} to show that
  $\tau\cdot\pi\in\O_I$ whenever $\pi\in\O_I$ and $\tau\in \hG$.  If
  $\pi\in\O_I$ then there exists an ideal $J$ of $A$ with continuous
  trace such that $\pi\in\O_J$. Note that
  $\tau\cdot\pi\in\tau\cdot\O_J=\O_{\tau\cdot J}$, where $\tau\cdot J
  = \hat\alpha_{\tau}(J)$.  Since $J$ has continuous trace each
  element $\rho$ of $\O_J$ is a Fell point and $\O_J$ is Hausdorff.
  Thus $\tau\cdot\O_J$ is also Hausdorff, and each point
  $\tau\cdot\rho$ in $\tau\cdot\O_J$ is a Fell point.  It follows that
  $\tau\cdot J$ is an ideal of $A$ with continuous trace and
  $\tau\cdot\pi\in\O_{\tau\cdot J}$.
  
  Finally, if $\set{ F_\alpha}$ is a neighborhood basis of $\pi$
  consisting of closed sets then $\set{ \tau\cdot F_\alpha}$ is a
  neighborhood basis of $\tau\cdot\pi$ with the same properties. Thus
  $\tau\cdot\pi\in\O_I$ by Lemma~\ref{lem-green}.
  
  We have shown that $\O_I$ and hence $I$ are $\hG$-invariant.  The
  final assertion follows from Lemma~\ref{lem-invar-ideals}.
\end{proof}

More generally, for an amenable $C^*$-dynamical system $(A,
G,\alpha)$, an ideal $I$ of $A\rtimes_\alpha G$ is invariant under the
dual coaction if and only if $I=J\rtimes_\alpha G$ for some unique,
$\alpha$-invariant ideal in $J$ of $A$
\cite[Theorem~3.4]{goot:_coaction}. Since we use a representation
theoretic approach to identify $\overline{\M(\cocx)}$ there are two
obstacles to extending our techniques to non-abelian groups. First,
there is no notion of induced coaction on $\speccocx$, and second, we
do not have a concrete description of $\speccocx$ in terms of $X$ and
$G$.

If $G$ is abelian, consider the quotient space obtained from
$X\times\hG$ where
\begin{equation*}
  (x,\omega)\sim(y,\tau) 
  \quad\text{if and only if}\quad\text{$\overline{G\cdot
  x}=\overline{G\cdot y}$ and $\omega\restr{S_x}=\tau\restr{S_y}$}.
\end{equation*}
This identification makes sense because $\overline{G\cdot
  x}=\overline{G\cdot y}$ implies $S_x=S_y$ for abelian groups.
Since we're assuming $(G,X)$ is second countable,
\cite[Theorem~5.3]{will81:_ccr} implies that
\[ 
[(x,\omega)]\mapsto\ker(\Ind_{(x,S_x)}^G(\omega\restr{S_x}))
\]
is a homeomorphism of $X\times\hG/\!\!\!\sim$ onto
$\Prim(C_0(X)\rtimes G)$.  We write $\pi_x(\omega)$ for
$\Ind_{(x,S_x)}^G(\omega\restr{S_x})$.  As noted in the paragraph
following the proof of \cite[Theorem~5.3]{will81:_ccr}, the map 
sending $(x,\omega)$ to $\ker \pi_{x}(\omega)$ is open from $X\times
\hG$ onto $\Prim\bigl(\cocx\bigr)$.  In particular, sets of the form
$U\times V/\!\!\sim$, with $U$ and $V$ open in $X$
and $\hG$, respectively, form a basis for the topology on
$\Prim\bigl(\cocx\bigr)$.

Let $\Sigma(G)$ denote the space of closed subgroups of $G$ endowed
with the compact Hausdorff topology from \cite{fell62:_topology}.  The
stability subgroups $S_x$ are said to \emph{vary continuously} if the
map $\sigma:X\to\Sigma(G): x\mapsto S_x$ is continuous.

If $A$ is a Fell algebra and $\pi\in\hat A$ then $\pi$ has an open
Hausdorff neighborhood in $\hat A$
\cite[Corollary~3.4]{ArchSom93:_trans}.  We want to be able to choose
this neighborhood to be $\hG$-invariant.

\begin{lemma}\label{lem-hausdorff} Suppose $(G,X)$ is a second
  countable transformation group with $G$ abelian and with
  continuously varying stability groups.  If $C_0(X)\rtimes G$ is a
  Fell algebra, then every irreducible representation of $\cocx$ has
  an open $\hG$-invariant Hausdorff neighborhood in $\speccocx$.
\end{lemma}

\begin{proof}
  Since $\cocx$ is postliminal, we can identify
  $\Prim\bigl(\cocx\bigr)$ and $\speccocx$.  We can view $\speccocx$
  as the appropriate quotient of $X/G\times\hG$, and then the map
  $(G\cdot x,\omega)\mapsto [\pi_{x}(\omega)]$ is an open surjection
  onto $\speccocx$ \cite[Theorem~5.3]{will81:_ccr}. In particular,
  (the class of) $\pi:=\pi_{x}(\omega)$ is a typical element of
  $\speccocx$. Since $A$ is a Fell algebra, $\pi$ has an open
  Hausdorff neighborhood \cite{ArchSom93:_trans} which is of the form
  $\O_J$ for some closed ideal $J$ of $A$.  We can shrink $J$ a bit if
  need be, and assume that there are open neighborhoods $U$ of $G\cdot
  x$ in $X/G$ and $V$ of $\omega$ in $\hG$ such that $U\times
  V/\!\!\sim$ is homeomorphic to $\O_J$. Suppose that $G\cdot x$ and
  $G\cdot y$ are distinct points in $U$.  Note that each orbit is
  closed in $X$ because $\cocx$ is liminal
  \cite[Theorem~3.1]{will81:_ccr}.  Thus, for each $\omega\in V$, the
  points $[G\cdot x,\omega]$ and $[G\cdot y,\omega]$ are distinct in
  $U\times V/\!\!\sim$.  Since $\O_{J}=U\times V/\!\!\sim$ is
  Hausdorff and $z\mapsto [G\cdot z,\omega]$ is continuous, we can
  separate $G\cdot x$ and $G\cdot y$ by $G$-invariant
  open sets and it follows that
  $U$ is Hausdorff.  Thus,
  \begin{equation*}
      \O:=U\times\hG/\!\!\sim
  \end{equation*}
  is a $\hG$-invariant neighborhood of $\pi$ which is Hausdorff
  because $U$ is Hausdorff and the stability subgroups vary
  continuously \cite{will82:_spectrum}.
\end{proof}

\section{Identifying ideals in $C_0(X)\rtimes G$}\label{sec:ident-ideals}

Let $(G,X)$ be a transformation group with continuously varying
stability groups. Define an equivalence relation on $X\times G$ by
\begin{equation*}
\text{$(x,s)\sim (y,t)$ if and only if $x=y$ and
$s^{-1}t\in S_x$.}
\end{equation*}
The continuity of the map $\sigma$ sending $x\mapsto S_x$ implies that
$X\times G/\!\!\sim$ is locally compact Hausdorff and that the
quotient map $\delta:X\times G\to X\times G/\!\!\sim$ is open
\cite[Lemma~2.3]{will81:_tr}.  The action of $G$ on $X$ is
\emph{$\sigma$-proper} if the map $[(x,s)]\mapsto (x,s\cdot x)$ of
$X\times G/\!\!\sim$ into $X\times X$ is proper
\cite[Definition~4.1]{rae88:_unitary_actions}.  It is not hard to see
that the action is $\sigma$-proper if and only if, given any compact
subset $K$ of $X$, the image in $X\times G/\!\!\sim$ of
\begin{equation}\label{eq:1}
\set{ (x,s)\in X\times G:\text{$x\in K$ and $s\cdot x\in K$}}
\end{equation}
is relatively compact. Any set $K$ for which the image of~\eqref{eq:1}
is relatively compact is called \emph{$G$-wandering}
\cite[p.~406]{rae88:_unitary_actions}.  If the action is free, then
the notions of $\sigma$-properness and $G$-wandering reduce to the
standard notions of properness and wandering, respectively.

\begin{lemma}
\label{lem-s-proper}
Let $(G,X)$ be a transformation group with continuously varying
stability groups.  If $U$ is an open $G$-wandering neighborhood of $X$
then the action of $G$ on $G\cdot U$ is $\sigma$-proper.
\end{lemma}

\begin{proof} 
  Let $K$ be a compact set in $G\cdot U$ and choose $t_1,\dots, t_n\in
  G$ such that $K\subset\bigcup_{i=1}^n t_i\cdot U$.  It suffices to
  show that for each $i$ and $j$,
\begin{equation}\label{eq:2}
\delta\bigl(
\set{(y,w)\in G\cdot U\times G: y\in K\cap t_i\cdot U \text{\ and\ }
w\cdot y\in K\cap t_j\cdot U}\bigr)
\end{equation}
is relatively compact in $(G\cdot U\times G)/\!\!\sim$.

Let $[(y_\alpha,w_\alpha)]$ be a net in the set described in
Equation~\ref{eq:2}. It will suffice to find a convergent subnet.
Since $\delta$ is open, we can pass to a subnet, relabel, and assume
that this net lifts to a net $(y_\alpha,s_\alpha)$ in $G\cdot U\times
G$ with $s_\alpha^{-1} w_\alpha\in S_{y_\alpha}$. Now $y_\alpha\in
K\cap t_i\cdot U$ and $s_\alpha\cdot y_\alpha\in K\cap t_j\cdot U$, so
that $y_\alpha= t_i\cdot x_\alpha$ for some $x_\alpha\in U$ and
$s_\alpha t_i\cdot x_\alpha=s_\alpha\cdot y_\alpha\in K\cap t_j\cdot
U$, that is, $t_j^{-1}s_\alpha t_i\cdot x_\alpha\in U$.

Now $\set{(x_\alpha,t_j^{-1}s_\alpha t_i)}$ is a net in $\set{ (y,w):
  y\in U \mbox{ and } w\cdot y\in U}$. Since $U$ is $G$-wandering
$\delta\bigl(\set{ (y,w): \text{$y\in U$ and $w\cdot y\in U$}}\bigr)$
is relatively compact.  By passing to a subnet and relabeling, we may
assume that for some $n_\alpha\in S_{x_\alpha}$ the net $\set{
  (x_\alpha, t_j^{-1}s_\alpha t_i n_\alpha)}$ converges in $X\times
G$. Since $t_{i}$ and $t_{j}$ are fixed,
\begin{equation*}
\set{ (t_i\cdot x_\alpha,
  s_\alpha t_i n_\alpha t_i^{-1})}=\set{( y_\alpha, s_\alpha t_i
  n_\alpha t_i^{-1})}
\end{equation*} 
also converges. Since $t_i n_\alpha t_i^{-1}\in S_{y_\alpha}$ and
$s_\alpha^{-1} w_\alpha\in S_{y_\alpha}$, we conclude that
$\set{[(y_\alpha, w_\alpha)]}$ converges.
\end{proof}

In \cite{aah:fell} the first author showed that if the action of $G$
on $X$ is free then $C_0(X)\rtimes G$ is a Fell algebra if and only if
$X$ is a Cartan $G$-space (that is, each point of $X$ has a wandering
neighborhood). If the stability subgroups vary continuously, we can
prove a similar result using the following generalization of
\cite[Proposition~1.1.4]{palais}.
\begin{lemma}
  \label{lem-closed-orbits}
  Suppose that $(G,X)$ is a second countable transformation group with
  $G$ abelian and with continuously varying stability groups.  If each
  point of $X$ has a $G$-wandering neighborhood, then $G\cdot x$ is
  closed in $X$ for all $x\in X$.
\end{lemma}
\begin{proof}
  Suppose that $y\in \overline{G\cdot x}$.  Let $U$ be a $G$-wandering
  neighborhood of $y$.  Then there are $s_{\alpha}\in G$ such that
  $s_{\alpha}\cdot x\to y$ and $s_{\alpha}\cdot x\in U$ for all
  $\alpha$.  We may replace $x$ by $s_{\alpha_{0}}\cdot x$ for some
  $s_{\alpha_{0}}\in G$, and assume that $x\in U$.  Then
  \begin{equation}
    \label{eq:3}
    \set{(x,s_{\alpha})}\subset \set{(z,s)\in X\times G: \text{$z\in
    U$ and $s\cdot z\in U$}}.
  \end{equation}
  Since the right-hand side of \eqref{eq:3} has relatively compact
  image in $X\times G/\!\!\sim$ and $\delta$ is open, we can pass to a
  subnet and relabel so that there are $t_{\alpha}\in S_{x}$ such that
  $s_{\alpha}t_{\alpha} \to s$ in $G$.  Then $y=s\cdot x$ and $G\cdot
  x$ is closed.
\end{proof}

\begin{prop}
\label{prop-fell}
Let $(G,X)$ be a second countable transformation group.  Suppose that
either $G$ acts freely, or that $G$ is abelian and that the stability
groups vary continuously. Then $C_0(X)\rtimes G$ is Fell algebra if
and only if each point of $X$ has a $G$-wandering neighborhood.
\end{prop}

\begin{proof}
  The free case is treated in \cite{aah:fell}.  Now suppose that $G$
  is abelian, that the stability groups vary continuously and that
  $\cocx$ is a Fell algebra. Fix $x\in X$ and let
  $\pi=\pi_{x}(1)\in\speccocx$.  By Lemma~\ref{lem-hausdorff}, $\pi$
  has an open Hausdorff $\hG$-invariant neighborhood $\O_J$, where $J$
  is an ideal of $A$.  Thus $J=C_0(Y)\rtimes G$ for some $G$-invariant
  open subset $Y$ of $X$, and $J$ has continuous trace.  The action of
  $G$ on $Y$ is $\sigma$-proper by
  \cite[Theorem~5.1]{will81:_tr}. Note that
  $x\in Y$, and let $N$ be a neighborhood of $y$ which is compact in
  $Y$. Then $N$ is $G$-wandering relative to $Y$, and since $Y$ is
  $G$-invariant $N$ is also $G$-wandering relative to $X$.
  
  Conversely, assume each point in $X$ has a $G$-wandering
  neighborhood.  Then Lemma~\ref{lem-closed-orbits} implies that the
  orbits are closed, and $\cocx$ is postliminal \cite{goot:_type}
  (even liminal \cite{will81:_ccr}).  In particular, each
  $\pi\in\speccocx$ is of the form $\pi_{x}(\omega)$ for some $x\in X$
  and $\omega\in\hG$.  Let $U$ be a $G$-wandering open neighborhood of
  $x$. By Lemma~\ref{lem-s-proper} the action of $G$ on $G\cdot U$ is
  $\sigma$-proper. Since the stability subgroups vary continuously it
  follows from \cite[Theorem~5.1]{will81:_tr} that $J=C_0(G\cdot
  U)\rtimes G$ is an ideal of $A$ which has continuous trace. Thus
  $\pi_x(\omega)$ is a Fell point of $\hat J$, whence it is also a
  Fell point of $\O_J\subset\hat A$.
\end{proof}

\begin{cor} 
  \label{cor-fell}
  Let $(G,X)$ be a second countable transformation group.  Suppose
  that either $G$ acts freely and $\cocx$ is EH-regular, or that $G$
  is abelian and that the stability groups vary continuously.  Then
  the largest Fell ideal of $\cocx$ is $C_0(W)\rtimes G$ where $W$ is
  the open $G$-invariant subset
  \begin{equation*}
     W=\set{ w\in X:\text{$w$ has a $G$-wandering neighborhood in $X$}}. 
  \end{equation*}
\end{cor}

\begin{proof}
  Again, the free case is dealt with in \cite{aah:fell}.  In any
  event, the largest Fell ideal of $\cocx$ is $J$ where
  $\O_J=\set{\pi\in\speccocx: \text{$\pi$ is a Fell point of
      $\speccocx$}}$.  Since $\O_J$ is invariant under the dual
  action, it follows that $J=C_0(W)\rtimes G$ for some open
  $G$-invariant subset $W$ of $X$.  Now apply
  Proposition~\ref{prop-fell}.
\end{proof}

\begin{remark}\label{rem-fu}
  When the action of $G$ on $X$ is free and $\cocx$ is postliminal,
  Green \cite[Corollary~18]{green77:_orbit} characterized the ideal
  $\overline{\M(\cocx)}$ as $C_{0}(Y')\rtimes G$ where
\begin{multline}\label{eq:10}
  Y'=\lbrace x\in X: \text{$x$ has a compact wandering neighborhood $N$}\\
  \text{such that $G\cdot N$ is closed in $X$}\,\rbrace;
\end{multline}
the following example shows that this is not quite correct. The
correct statement is contained in Theorem~\ref{thm-new-main} below and
says that the open subset $Y$ of $X$ corresponding to
$\overline{\M(\cocx)}$ is given by Equation~\ref{eq:7} in \S1.
\end{remark}

\begin{example}\label{ex-palais}
  Consider the transformation group described by Palais in
  \cite[p.~298]{palais}, where $X$ is the strip
  $\set{(x,y):\text{$-1\leq x\leq 1$ and $ y\in\R$}}$ and the group
  action is by $G=\R$.  Beyond the strip $-1<x<1$ the action moves a
  point according to
  \begin{equation*}
  t\cdot(1,y)=(1,y+t)\quad\text{and}\quad t\cdot(-1,y)=(-1, y-t).
  \end{equation*}
  If $(x_0,y_0)\in\interior(X)$ let $C_{(x_0,y_0)}$ be the vertical
  translate of the graph of $y=\frac{x^2}{1-x^2}$ which passes through
  $(x_0,y_0)$.  Define $t\cdot(x_0,y_0)$ to be the point $(x,y)$ on
  $C_{(x_0,y_0)}$ such that the length of the arc of $C_{(x_0,y_0)}$
  between $(x_0,y_0)$ and $(x,y)$ is $|t|$, and $x-x_0$ has the same
  sign as $t$.  That is, $(x_0,y_0)$ moves counter-clockwise along
  $C_{(x_0,y_0)}$ at unit speed.

Palais states that a compact set is wandering if and only if it meets
at most one of the lines $x=1$ and $x=-1$; this is only partially
correct. Certainly, if a compact set meets at most one of the boundary
lines then it is wandering. However, $N=[0,1]\times [-1,1]\cup\set{
  (-1,0)}$ is an example of a wandering compact set meeting both
boundary lines; moreover, $G\cdot N$ is closed in $X$, and $N$ is a
neighborhood of $(1,y)$ for all $y\in(-1,1)$.  One sees
from these examples that for this transformation group, the set $Y'$
described in \eqref{eq:10} is all of $X$ whence $C_0(X)\rtimes G$
should have continuous trace. But this is impossible because
$X/G\cong\speccocx$ is not Hausdorff: for example, $G\cdot(-1+1/n,0)$
is a sequence which converges to the distinct orbits $G\cdot(-1,0)$
and $G\cdot(1,0)$. Alternatively, note that not every compact set is
wandering which contradicts \cite[Theorem~17]{green77:_orbit}.
\end{example}

\begin{remark}
  \label{rem-wandering}
  In Theorem~\ref{thm-new-main}, we want to consider sets $K\subset X$
  which are $G$-wandering even though we definitely are not assuming
  that the stabilizer map $\sigma$ is continuous on all of $X$.  To
  make sense of this, we have to assume that $\sigma$ is at least
  continuous on
  $G\cdot K$, and then it makes sense to ask if $K$ is $G$-wandering in
  $G\cdot K$ (or, equivalently, in any $G$-invariant set $Z$ which
  contains $K$ and on which $\sigma$ is continuous). If $K$ is open,
  it is not hard to see that $\sigma$ is continuous on $G\cdot K$ if
  and only if $\sigma$ is continuous on $K$.  However, in general the
  continuity of $\sigma$ on $K$ does not imply that $\sigma$ is
  continuous on $G\cdot K$.  The next lemma will allow us to ignore
  this difficulty when applying the Theorem.
\end{remark}

\begin{lemma}
  \label{lem-wandering}
  Suppose that $(G,X)$ is a locally compact
  transformation group with $G$ abelian and with stabilizer map
  $\sigma$.  Let $q:X\to X/G$ be the quotient map.  If $\sigma$ is
  continuous on a compact set $K$ and if $q(K)$ is Hausdorff, then
  $\sigma$ is continuous on $G\cdot K$.
\end{lemma}
\begin{proof}
  Suppose that $r_{\alpha}\cdot x_{\alpha}\to r\cdot x$ for
  $r_{\alpha},r\in G$ and $x_{\alpha},x\in K$.  We want to show that
  $S_{r_{\alpha}\cdot x_{\alpha}}=S_{x_{\alpha}}$ converges to
  $S_{r\cdot x}=S_{x}$.  Since this happens if and only if every
  subnet converges to $S_{x}$, we can pass to a subnet, relabel and
  assume that $x_{\alpha}\to y\in K$.  Since $q(K)$ is Hausdorff,
  $y=s\cdot x$ for some $s\in G$.  Thus by assumption,
  $S_{r_{\alpha}\cdot x_{\alpha}}=S_{x_{\alpha}}$ converges to
  $S_{y}=S_{x}$. 
\end{proof}

\begin{thm}\label{thm-new-main} 
  Let $(G,X)$ be a second countable transformation group, and let
  $\sigma$ be the stabilizer map sending $x\mapsto S_{x}$.  Assume
  either that $G$ acts freely and $\cocx$ is EH-regular, or that $G$
  is abelian.  Let $I:=\overline{\M\bigl(\cocx\bigr)}$.  Then
  $I=C_0(Y)\rtimes G$, where $Y$ is the open $G$-invariant subset
  \begin{multline}\label{eq:4}
    Y=\lbrace\, y\in X: \text{$\sigma$ is continuous on a $G$-wandering compact}\\
    \text{neighborhood $N$ of $y$ such that $q(N)$ is closed and
      Hausdorff}\,\rbrace,
  \end{multline}
  where $q:X\to X/G$ is the quotient map.
 \end{thm}

\begin{proof} 
  Our proof is modeled on the proof of
  \cite[Corollary~18]{green77:_orbit}.  Here we'll give the proof for
  $G$ abelian and remark that the free case follows from the same sort
  of argument together with the following observation. If the action
  is free, then $EH$-regularity implies that
  $\speccocx$ is homeomorphic to the
  $\mathrm{T}_{0}$-ization $(X/G)^{\sim}$ of $X/G$
  \cite[Corollary~5.10]{will81:_ccr}. It follows that the map
  $Y\mapsto \cocy$ from the set of $G$-invariant open subsets of $X$
  to the set of ideals of $\cocx$ is a bijection.
  
  By Proposition~\ref{prop-cts-trace-id}, $I=C_0(Z)\rtimes G$ where
  $Z$ is an open $G$-invariant subset of $X$. Let $Y$ be as in
  \eqref{eq:4}.  Suppose that $\pi\in \O_I$.  Since $I$ has continuous
  trace, it is certainly postliminal, and $\pi=\pi_{x}(\omega)$ for
  $x\in Z$ and $\omega\in\hG$. Furthermore,
  \cite[Theorem~5.1]{will81:_tr} implies that the stabilizer map
  $\sigma$ is continuous on $Z$ and that the action of $G$ on $Z$ is
  $\sigma$-proper. Let $N$ be a compact neighborhood of $x$ in $Z$.
  Then $N$ is $G$-wandering relative to $Z$, and since $Z$ is
  $G$-invariant, $N$ is also $G$-wandering relative to $X$.
  
  Let $q:X\to X/G$ be the quotient map.  We claim there is a closed
  neighborhood $V$ of $G\cdot x$ in $X/G$ such that $V\subset q(N)$.
  To prove the claim, we identify $\Prim\bigl(\cocx\bigr)$ with
  $X\times \hG/\!\!\sim$.  Then Lemma~\ref{lem-green} implies
  $\ker\pi_{x}(\omega)$ has a closed neighborhood $W\subset
  (N\times\hG)/\!\!\sim$.  The map $y\mapsto \ker\pi_{y}(\omega)$ is
  continuous by \cite[Lemma~4.9]{will81:_ccr}, and factors through
  $X/G$ by \cite[Corollary~4.8]{will81:_ccr}. Thus we get a continuous
  map $s_{\omega}: X/G\to \Prim\bigl(\cocx\bigr)$.  Let
  $V:=s_{\omega}^{-1}(W)$.  Then $V$ is a closed neighborhood of
  $G\cdot x$.  To prove the claim, it remains to see that $V\subset
  q(N)$.  But if $G\cdot y\in V$, then there is a $(z,\gamma)\in
  N\times\hG$ such that $(y,\omega)\sim(z,\gamma)$.  In particular,
  $\overline{G\cdot y}=\overline{G\cdot z}$.  Since $Z$ is open and
  $G$-invariant, it follows that $y\in Z$. (We have $s_{\alpha}\cdot y
  \to z$ for $s_{\alpha}\in G$.)  Thus $G\cdot y$ and $G\cdot z$ have
  the same closures in $Z$.  But $C_{0}(Z)\rtimes G$ is liminal and
  each orbit must be closed in $Z$
  \cite[Proposition~4.17]{will81:_ccr}.  Thus $G\cdot y=G\cdot z\in
  q(N)$ as claimed.
  
  With $V$ as above, set $N'=q^{-1}(V)\cap N$.  Note that $N'$
  is compact and $G$-wandering and $G\cdot N'=q^{-1}(V)$ is closed.
  Finally, $G\cdot N'/G$ is Hausdorff because $G\cdot N'\subset Z$,
  and $Z/G$ is Hausdorff since $C_0(Z)\rtimes G$ has continuous trace
  \cite{will82:_spectrum}.  This implies that $x\in Y$.  Therefore
  $Z\subset Y$, and $I=C_{0}(Z)\rtimes G\subset C_{0}(Y)\rtimes G$.
  
  To prove the reverse implication notice that $C_{0}(Y)\rtimes G$ is
  a Fell algebra by Proposition~\ref{prop-fell}.  In particular, it is
  postliminal, and every irreducible representation of
  $C_{0}(Y)\rtimes G$ is of the form $\pi=\pi_{y}(\omega)$ for $y\in
  Y$ and $\omega\in\hG$.  We will show that $\pi\in\O_I$ by verifying
  items~(\ref{item:1}) and (\ref{item:2}) of Lemma~\ref{lem-green}.
  Since $C_{0}(Y)\rtimes G$ is a Fell algebra $\pi$ has a Hausdorff
  open neighborhood $\O_J$, where $J$ is a closed ideal of
  $C_{0}(X)\rtimes G$ \cite[Corollary~3.4]{ArchSom93:_trans}.  Note
  that $J$ is a Fell algebra with Hausdorff spectrum.  Hence $J$ has
  continuous trace.  This establishes item~(\ref{item:1}) of
  Lemma~\ref{lem-green}.
  
  Let $N$ be a compact $G$-wandering neighborhood of $y$ as
  in~\eqref{eq:4}.  We identify $\speccocx$ with
  $X\times\hG/\!\!\sim$.  Note that $V=G\cdot N\times \hG/\!\!\sim$ is
  a closed neighborhood of $\pi$ (first consider the complement and
  recall that the quotient map is open).  That $V$ is Hausdorff
  follows from \cite{will82:_spectrum} because $G\cdot N/G$ is
  Hausdorff and the stability subgroups vary continuously on $G\cdot
  N$ by Lemma~\ref{lem-wandering}.  Let $\set{ F_\alpha}$ be a
  neighborhood basis of $\pi$ in $\speccocx$ consisting of compact
  sets.  Since a compact subset of a Hausdorff space is closed, $\set{
    F_\alpha\cap V}$ is a neighborhood basis of $\pi$ in $\speccocx$
  consisting of closed sets.  This establishes item~(\ref{item:2}).
  Since $\pi$ was an arbitrary irreducible representation of
  $C_{0}(Y)\rtimes G$, we must have $C_{0}(Y)\rtimes G\subset
  I=C_{0}(Z)\rtimes G$.  Therefore $Z=Y$ and we're done.
\end{proof}

\begin{example}
  \label{ex-palias-crt}
  If $A=C_{0}(X)\rtimes\R$ is the transformation group in
  Example~\ref{ex-palais}, then $I=\overline{\M(A)}$ corresponds to
  the open strip $Y=\set{(x,y):-1<x<1}$.
\end{example}

\begin{example}
  Let $G=\R^+$ act on $X=\R^2$ by $t\cdot (x,y)=(x/t, y/t)$. The
  orbits are rays emanating from the origin together with the origin
  which is a fixed point. Each orbit is locally closed so $\cocx$ is
  postliminal \cite{goot:_type}.  The stability subgroups do not vary
  continuously on any neighborhood of $(0,0)$. If $U$ is any
  $G$-wandering (hence wandering) neighborhood of $(x,y)\neq(0,0)$
  then $(0,0)\in \overline{G\cdot U}$ so that $G\cdot U$ is not closed
  in $X$. Thus Theorem~\ref{thm-new-main} implies that
  $\M(\cocx)=\set0$. Note that the action of $G$ on
  $W:=X\smallsetminus\set{ (0,0)}$ is free and proper so that
  $C_{0}(W) \rtimes G$
  is an essential ideal of $\cocx$ with continuous trace.
\end{example}

It should be pointed out that even for liminal algebras $A$, it is
possible that $\M(A)=\set0$.  To see this, recall that a
point $x$ of a topological space $X$ is \emph{separated} if for any
point $y$ of $X$ not in the closure of $\set{x}$, the points $x$ and
$y$ admit a pair of disjoint neighborhoods.  If $A$ is a separable
$C^*$-algebra, then the set $\calS$ of separated points of the
spectrum $\hat A$ is a dense $G_\delta$ \cite[3.9.4]{dix69:_book}.

\begin{lemma}\label{lem-containment}
  Let $A$ be a $C^*$-algebra and $I:=\overline{\M(A)}$.  Then $\O_I$ is
  contained in the interior of the separated points $\calS$ of $\hat
  A$.
\end{lemma}

\begin{proof}
  Let $\pi\in\O_I$, and $\rho\in\hat A$ such that
  $\rho\notin\overline{\set{\pi}}$.
  If $\rho\in\O_I$ then $\rho$ and $\pi$ can be separated by disjoint
  relative open subsets of $\hat A$ because $\O_I$ is Hausdorff. Since
  $\O_I$ is open these relative open sets are open.  Now suppose that
  $\rho\notin\O_I$. Fix a positive element $a$ of $\M(A)$ such that
  $\tr(\pi(a))>1$ and let $f:\hat A\to [0,\infty)$ be the (continuous)
  map $\sigma\mapsto\tr(\sigma(a))$. Note that $\rho(a)=0$. Now
  $f^{-1}((1,\infty))$ and $f^{-1}([0,\frac{1}{2}))$ are disjoint open
  neighborhoods of $\pi$ and $\rho$, respectively.  Thus
  $\O_I\subset\calS$ and since $\O_I$ is open we have
  $\O_I\subset\interior\calS$.
\end{proof}

Dixmier has given an example of a separable liminal $C^{*}$-algebra
$A$ such that the interior of the separated points in $\hat A$ is empty
\cite[Proposition 4]{dix61:_separated}. 
Thus $\M(A)=\set0$ for this algebra.

\begin{thm}\label{thm-liminal-ideal}
  Let $(G,X)$ be a second countable transformation group.  Suppose
  that either $G$ acts freely and $\cocx$ is EH-regular, or that $G$
  is abelian.  Then the largest liminal ideal of $C_0(X)\rtimes G$ is
  $C_0(Z)\rtimes G$ where $Z$ is the open $G$-invariant subset
\begin{multline}\label{eq-ccr}
  Z=\lbrace\, x\in X: \text{$x$ has a neighborhood $U$}\\ \text{such
    that $G\cdot z$ is closed in $G\cdot U$ for each $z\in
    U$}\,\rbrace.
\end{multline}
\end{thm}
 
\begin{proof}  If $J$ is the largest liminal ideal then
  $\O_J=\set{\pi\in\hat A:\pi\bigl(\cocx\bigr)=\K(\H_\pi)}$.  If $G$
  is abelian then $\O_{J}$ is invariant under the dual action, and we
  have $J=C_0(Y)\rtimes G$ for some open $G$-invariant subset $Y$ of
  $X$.  This is trivial in the free case.  Let $Z$ be as in~\eqref{eq-ccr}. Note that every $y\in Y$
  has a neighborhood $U$ (namely $Y$) such that $G\cdot z$ is closed
  in $G\cdot U$ for every $z\in U$ by \cite[Theorem~3.1]{will81:_ccr}, so $Y\subset Z$.
  
  Let $x\in Z\smallsetminus Y$. Let $V$ be an open
  neighborhood of $x$ such that $G\cdot z$ is closed in $G\cdot V$ for
  each $z\in V$.  Not every orbit in $Y'=Y\cup G\cdot V$ can be closed
  in $Y'$ because $C_0(Y)\rtimes G$ is the largest liminal ideal. Suppose that $G\cdot
  z$ is not closed in $Y'$.  Then there exists $s_\alpha\in G$ and
  $w\in Y'$ such that $s_\alpha\cdot z\to w\notin G\cdot z$.
  
  Since $w\in Y'$, $w$ has a neighborhood $W$ such that $G\cdot u$ is
  closed in $G\cdot W$ for all $u\in W$.  But we can assume that
  $s_{\alpha_{0}}\cdot z\in W$ for some $s_{\alpha_{0}}$ and then
  $G\cdot s_{\alpha_{0}}\cdot z=G\cdot z$ must be closed in $G\cdot
  W$.  Thus $w\in G\cdot z$, and this is a contradiction. Hence $Z=Y$ and we are done.
\end{proof}

Every $C^*$-algebra $A$ has a largest postliminal ideal $I$, and this
ideal $I$ is the smallest ideal such that the corresponding quotient
is anti-liminal \cite[Proposition~4.3.6]{dix69:_book}. When
$A=C_0(X)\rtimes G$ and $G$ is abelian, it is clear that $I$ is
invariant under the dual action: for every $\tau\in\hG$ the ideal
$\hat\alpha_\tau(I)$ is postliminal and $A/{\hat\alpha_\tau(I)}$ is
antiliminal, hence $\hat\alpha_\tau(I)\subset I$. If $G$ is abelian or
$G$ acts freely then $C_0(X)\rtimes G$ is Type I if and only if $X/G$
is $\mathrm{T}_0$ \cite[Theorem~3.3]{goot:_type}.  Effros and Glimm
have given a number of conditions on a second countable locally
compact transformation group $(G,X)$ which are equivalent to $X/G$
being $\mathrm{T}_{0}$: see \cite{gl1},
\cite[Theorems~2.1~and~2.6]{effros} and \cite{effros-gpoid}.  For
example, $X/G$ is $\mathrm{T}_{0}$ if and only if each orbit is
\emph{regular}: the map $sS_{x}\mapsto s\cdot x$ is a homeomorphism of
$G/S_{x}$ onto $G\cdot x$.\footnote{The term regular is borrowed from
  the Definition on page~223 of \cite{green78}.}  Using the
Effros-Glimm results, we have the following.

\begin{lemma} [Glimm-Effros] \label{lem-glimm-effros}
  Suppose that $(G,X)$ is a second countable locally compact
  transformation group and that $U$ is a neighborhood of $x\in X$.
  Then the following are equivalent.
  \begin{enumerate}
  \item $G\cdot U/G$ is $\mathrm{T}_{0}$ in the quotient topology.
  \item $G\cdot y$ is regular for each $y\in U$.
  \item $G\cdot y$ is a $G_{\delta}$ subset of $X$ for each $y\in U$.
  \item $G\cdot y$ is locally closed in $X$ for each $y\in U$.
  \item $G\cdot y$ is second category in itself for each $y\in U$.
  \end{enumerate}
\end{lemma}

\begin{thm}\label{thm-I-ideal}
  Let $(G,X)$ be a second countable transformation group. Suppose that
  either $G$ acts freely and $\cocx$ is EH-regular, or that $G$ is
  abelian.  Then the largest postliminal ideal of $C_0(X)\rtimes
  G$ equals $C_0(Z)\rtimes G$ where $Z$ is the $G$-invariant subset
\begin{equation}\label{eq-gcr}
Z=\set{ x\in X:\text{$ x$ has a neighborhood $U$ such that
$G\cdot U/G$ is $\mathrm{T}_0$}}.
\end{equation}
\end{thm}
\begin{remark}
  The set $Z$ can be realized as the set of points with neighborhoods
  satisfying any of the equivalent conditions of
  Lemma~\ref{lem-glimm-effros}.
\end{remark}

\begin{proof}
  If $G$ is abelian, the largest postliminal ideal of
  $C_0(X)\rtimes G$ is invariant under the dual action, so equals
  $C_0(Y)\rtimes G$ for some $G$-invariant open subset $Y$ of $X$.
Let $Z$ be as in \eqref{eq-gcr}.
  Every $y\in Y$ has an open $G$-invariant neighborhood $U$ (namely
  $Y$) such that $G\cdot U/G$ is ${\rm T}_0$ by
  \cite[Theorem~3.3]{goot:_type}. Thus $Y\subset Z$.
  
  Let $x\in Z\smallsetminus Y$ and  $V$  an open
  neighborhood of $x$ such that $G\cdot V/G$ is ${\rm T}_0$. Note that
  $T:=(G\cdot V\cup Y)/G$ cannot be ${\rm T}_0$ by the maximality of
  $C_0(Y)\rtimes G$. Choose distinct points $G\cdot z_1$ and $G\cdot
  z_2$ in $T$ such that every open neighborhood $U_1$ of $G\cdot z_1$
  contains $G\cdot z_2$ and every open neighborhood $U_2$ of $G\cdot
  z_2$ contains $G\cdot z_1$.
  
  If $G\cdot z_1\in T\smallsetminus (G\cdot V/G)$ and $G\cdot z_2\in
  T\smallsetminus (Y/G)$ then $G\cdot V/G$ is an open neighborhood of
  $G\cdot z_2$ which does not contain $G\cdot z_1$, which is a
  contradiction.
  
  If $G\cdot z_{1}$ and $G\cdot z_{2}$ both belong to $Y/G$
  or if $G\cdot z_{1}$ and $G\cdot z_{2}$ both belong to $G\cdot V/G$
  then we get an immediate contradiction because $Y/G$ and $G\cdot
  V/G$ are open and ${\rm T}_0$. Hence $Y=Z$.
\end{proof}

\begin{remark}
  Let $(A, G, \alpha)$ be a $C^*$-dynamical system with $G$ compact
  (but not necessarily abelian). It follows from
  \cite[Propositions~2.3~and~2.5]{goot:_compact} that the largest
  liminal and postliminal ideals in $A\rtimes_\alpha G$ are of the
  form $J\rtimes_\alpha G$ where $J$ is an $\alpha$-invariant ideal of
  $A$. This is trivial if $A=C_0(X)$, because $\cocx$ has bounded
  trace (hence is liminal) when $G$ is compact
  \cite[Proposition~3.4]{aah:trace}.
\end{remark}

\bibliographystyle{amsplain}


\providecommand{\bysame}{\leavevmode\hbox to3em{\hrulefill}\thinspace}

\end{document}